\def\k{{\kappa}}
\def\Th{{\Theta}}
\def\w{{\omega}}
\def\W{{\Omega}}

\def\p{{\bm p}}
\def\N{{\bm N}}
\def\Z{{\bm Z}}
\def\bb{{\overline b}}
\def\c{{\overline c}}
\def\d{{\overline d}}
\def\e{{\overline e}}
\def\ep{{\epsilon}}
\def\u{{\overline u}}
\def\vv{{\overline v}}
\def\0{{\overline 0}}

\def\C{{\cal C}}
\def\F{{\cal F}}
\def\G{{\cal G}}
\def\cH{{\cal H}}
\def\K{{\cal K}}
\def\M{{\cal M}}
\def\P{{\cal P}}
\def\Q{{\cal Q}}
\def\S{{\cal S}}

\def\Wh{{\cal W}}

\def\inf{{\infty}}
\def\mt{{\emptyset}}
\def\rar{{\rightarrow}}
\def\la{{\langle}}
\def\ra{{\rangle}}

\def\lf{{\lfloor}}
\def\rf{{\rfloor}}
\def\pr{{^\prime}}
\def\sn{{n^{1/2}}}

\documentclass[12pt]{article}
\usepackage{latexsym,psfig,multicol}

\newcommand{\bm}[1]{{\mbox{\boldmath $#1$}}}
\newtheorem{theorem}{Theorem}[section]
\newtheorem{conjecture}[theorem]{Conjecture}
\newtheorem{corollary}[theorem]{Corollary}
\newtheorem{fact}[theorem]{Fact}
\newtheorem{lemma}[theorem]{Lemma}
\newtheorem{question}[theorem]{Question}

\begin{document}

% #######################################################################
% #######################################################################
%       
%       TITLE PAGE: 
%       
\title{A Survey of Graph Pebbling}
\author{
Glenn H. Hurlbert\\
Department of Mathematics\\
Arizona State University\\
Tempe, Arizona 85287-1804\\
email: hurlbert@math.la.asu.edu
}
\date{}
\maketitle
%        
%\newpage

\small\normalsize

% #######################################################################
% #######################################################################
%       
%       ABSTRACT:
%       
\begin{abstract}
We survey results on the pebbling numbers of graphs as well as their 
historical connection with a number-theoretic question of Erd\H os and Lemke.
We also present new results on two probabilistic pebbling considerations, 
first the random graph threshold for the property that the pebbling number 
of a graph equals its number of vertices, and second the pebbling threshold 
function for various natural graph sequences.
Finally, we relate the question of the existence of pebbling thresholds to 
a strengthening of the normal property of posets, and show that the multiset 
lattice is not supernormal.
\vspace{0.2 in}

\vspace{0.2 in}

\end{abstract}
% 
%\small\normalsize
%\renewcommand{\baselinestretch}{1.00}
\setlength{\topmargin}{.25in}
\setlength{\textheight}{7.5in}
\newpage

% #######################################################################
% #######################################################################
%       
%       BEGINNING OF PAPER:
%       
\section{Introduction}\label{Intro}
Suppose $t$ pebbles are distributed onto the vertices of a graph $G$.
A pebbling step $[u,v]$ consists of removing two pebbles from one vertex $u$ 
and then placing one pebble at an adjacent vertex $v$.
We say a pebble can be {\it moved} to a vertex $r$, the {\it root} vertex,
if we can repeatedly apply pebbling steps so that in the resulting
distribution $r$ has one pebble.
The pebbling step $[u,v]$ is {\it greedy} if $dist(v,r)<dist(u,r)$, and
{\it semigreedy} if $dist(v,r)\le dist(u,r)$.
Here the function $dist$ denotes distance.

For a graph $G$, we define the {\it pebbling number}, $f(G)$, to be the
smallest integer $t$ such that, for any distribution of $t$ pebbles onto the
vertices of $G$, one pebble can be moved to any specified root vertex $r$.
If $D$ is a distribution of pebbles onto the vertices of $G$ and it is
possible to move a pebble to the root vertex $r$, then we say that
$D$ is $r$-{\it solvable}; otherwise, $D$ is $r$-{\it unsolvable}.
Then $D$ is {\it solvable} if it is $r$-solvable for all $r$, and
{\it unsolvable} otherwise.
We denote by $D(v)$ the number of pebbles on vertex $v$ in $D$ and let the
{\it size}, $|D|$, of $D$ be the total number of pebbles in $D$; that is,
$|D|=\sum_vD(v)$.
This yields another way to define $f(G)$, as one more than the maximum $t$ 
such that there exists an unsolvable pebbling distribution of size $t$.

Throughout this paper $G$ will denote a simple connected graph, where
$n(G)=|V(G)|$, and $f(G)$ will denote the pebbling number of $G$.
For any two graphs $G_1$ and $G_2$, we define the {\it cartesian product}
$G_1 \Box G_2$ to be the graph with vertex set
$V(G_1\Box G_2)=\{(v_1,v_2) | v_1 \in V(G_1), v_2 \in V(G_2)\}$ and edge set
$E(G_1\Box G_2)=\{((v_1,v_2),(w_1,w_2)) | (v_1=w_1$ and $(v_2,w_2)\in E(G_2))$
or $(v_2=w_2$ and $(v_1,w_1)\in E(G_1))$\}.
Thus the $m$-dimensional cube $Q^m$ can be written as the cartesian product
of $K_2$ with itself $m$ times.

In this paper we survey the current knowledge regarding the many questions
surrounding this simple pebbling operation.
Our aim is to tie together several of Paul Erd\H os's favorite subjects,
namely, graph theory, number theory, probability, and extremal set theory.
It is interesting that graph pebbling arose out of attempts to answer a
question of Erd\H os.

We begin in Section \ref{Pebbling} by reviewing the obvious general upper and 
lower bounds, the known results for the pebbling numbers of various classes 
of graphs, and conjectures and theorems involving products, diameter, and 
connectivity.
Section \ref{Number} details the origins of graph pebbling and its 
connections with number theory.
In Section \ref{Thresholds} we define threshold functions for graph sequences,
analogous to threshold functions for random graph properties, and discuss
known results and open problems.
Finally, Section \ref{Sets} addresses the underlying set theory on which 
the existence of pebbling threshold functions relies, and contrasts this
theory with the set theory underpinning the existence of random graph
threshold functions.

% #######################################################################
% #######################################################################
%       
%       NEW SECTION:
%       
\section{Pebbling Numbers}\label{Pebbling}
If one pebble is placed at each vertex other than the root vertex, $r$,
then no pebble can be moved to $r$.
Also, if $w$ is at distance $l$ from $r$, $2^{l}-1$ pebbles are placed at $w$, 
and no pebbles are placed elsewhere, then no pebble can be moved to $r$.
On the other hand, if more than $(2^d-1)(n-1)$ pebbles are placed on the
vertices of a graph of diameter $d$, then either every vertex has at least
one pebble on it or some vertex $w$ has at least $2^d$ pebbles on it.
In either case one can immediately pebble from $w$ to any vertex $r$.
We record these observations as

\begin{fact}\label{LowBd}
Let $d=diam(G)$ and $n=n(G)$.
Then $max\{n, 2^d\}\le f(G)\le (2^d-1)(n-1)+1$.
\end{fact}

Of course this means that $f(K_n)=n$, where $K_n$ is the complete graph on
$n$ vertices.
Let $P_n$ denote the path on $n+1$ vertices.
Here is a simple weight function argument to show that $f(P_n)=2^n$.
For a given distribution $D$ and leaf root $r$ define the weight 
$w(D)=\sum_vw(v)$, where $w(v)=D(v)/2^{dist(v,r)}$.
Because the weight of a distribution is preserved under greedy pebbling steps,
$D$ is an $r$-unsolvable distribution if and only if $w(D)<1$.
Because pebbling reduces the size of a distribution, if $D$ has maximum size 
with respect to $r$-unsolvable distributions then all its pebbles lie on the 
leaf opposite from $r$, implying $|D|=2^n-1$.
Finally, for any other choice of root $r$,
one applies the above argument to both neighbors of $r$ and
notices that $(2^a-1)+(2^b-1)<2^{a+b}-1$.

Let $C_n$ be the cycle on $n$ vertices.
It is easy to see that $f(C_5)=5$ and $f(C_6)=8$, so that each of the two 
lower bounds are relevant. 
The pebbling numbers of cycles is derived in \cite{Snevily}. 
In the case of larger odd cycles, the pebbling number exceeds
both lower bounds.

\begin{theorem}\cite{Snevily}\label{Cycles}
For $k\ge 1$, $f(C_{2k}) = 2^k$ and 
$f(C_{2k+1}) = 2\lfloor \frac{2^{k+1}}{3} \rfloor +1$.
\end{theorem} 

We can use this result to prove that the pebbling number of the 
Petersen graph $P$ is 10.
A common drawing of $P$ displays inner and outer 5-cycles.
Consider a distribution $D$ of size 10 with root $r$, having no pebble on it.
If a neighbor $s$ of $r$ has a pebble on it then by symmetry we may draw $P$
so that $r$ is on the outer cycle and $s$ is on the inner.
Because $f(C_5)=5$ we may assume that there are fewer than 5 pebbles on the
outer cycle, and thus more than 5 on the inner.
We may ignore one of the pebbles on $s$ and use 5 of the other pebbles to 
move a second pebble to $s$, then one to $r$.
On the other hand we consider the case that the neighbors of $r$ are also 
void of pebbles, so all 10 pebbles are on the 6-cycle formed by the 
nonneighbors of $r$.
From there one can see that 2 pebbles can be moved to a neighbor of 
$r$, then one to $r$.

The pebbling number of a tree $T$ on $n$ vertices is more complicated.
Consider a partition $Q=(Q_1,\ldots,Q_m)$ of the edges of $T$ into paths 
$Q_1,\ldots,Q_m$, written so that $q_i\ge q_{i+1}$, where $q_i=|Q_i|$.
Any choice of root vertex $r$ in $T$ induces an orientation of the edges 
of $T$ and thus also on each path $Q_i$.
The orientation on $Q_i$ determines a root $r_i$ of $Q_i$, which may or 
may not be an endpoint of $Q_i$.
If $r_i$ is an endpoint of $Q_i$ then we say that $Q_i$ is 
{\it well $r$-directed}.
We call $Q$ an {\it $r$-path partition} of $T$ if each path $Q_i$ is 
well $r$-directed, and a {\it path partition} if it is an $r$-path partition 
for some $r$.
The path partition $Q$ {\it majorizes} another, $Q\pr$, if its sequence of
path lengths majorizes that of the other, that is, if $q_j>q_j\pr$, where 
$j=\min\{i:q_i\ne q_i\pr\}$.
A path (resp. $r$-path) partition of $T$ is {\it maximum} (resp. 
{\it $r$-maximum}) if no other path (resp. $r$-path) partition majorizes it.

\begin{theorem}\cite{Moews}\label{Trees}
Let $(q_1,q_2,\ldots,q_m)$ be the nonincreasing sequence of path lengths of
a maximum path partition $Q=(Q_1,\ldots,Q_m)$ of a tree $T$.
Then $f(T)=\Biggl(\sum\limits_{i=1}^m 2^{q_i}\Biggr) -m+1$.
\end{theorem} 

The crucial idea in the argument is to find the right generalization to use
for an induction proof.
Say that a distribution is {\it $k$-fold $r$-solvable} if it is possible to
move $k$ pebbles to the vertex $r$ after a sequence of pebbling steps.
Define $f(G,r;k)$ to be the minimum $t$ so that every distribution of size
at least $t$ is $k$-fold $r$-solvable.
Moews \cite{Moews} proves that 
$f(T,r;k)=k2^{q_1}+\Bigl(\ \sum\limits_{i=2}^m 2^{q_i}\Bigr) -m+1$, 
where $(q_1,q_2,\ldots,q_m)$ is the nonincreasing sequence of path lengths 
of a maximum $r$-path partition $Q=(Q_1,\ldots,Q_m)$ of the tree $T$. 
If $T-r$ is the union of trees $T_1,\ldots,T_s$, with each $T_j$ rooted at
a neighbor $r_i$ of $r$, then he uses as the inductive hypothesis the equality
$$f(T,r;k) = \max\Biggl\{\ \sum\limits_{i=1}^s f(T_i,r_i;k_i+1)\Biggr\},$$
where the maximum is over all $k_1,\ldots,k_s$ which satisfy
$\sum\limits_{i=1}^s \lfloor k_i/2\rfloor <k$.

It is natural to ask what pebbling solutions can look like.
We say a graph $G$ is {\it greedy} ({\it semi-greedy}) if every distribution
of size at least $f(G)$ has a solution in which every pebbling step is greedy
(semi-greedy).
We say a graph $G$ is {\it tree-solvable} if every distribution of size at 
least $f(G)$ has a solution in which the edges traversed by pebbling steps 
form an acyclic subgraph.
The 5-cycle $abcde$ is not greedy, as witnessed by the distribution 
$D(c,d)=(3,2)$, with root $a$.
Worse yet, let $H$ be the graph formed from the 6-cycle $abcdef$ by adjoining
new vertices $g$ to $a$ and $c$, and $h$ to $a$ and $e$.
It is not difficult to show that $f(H)=9$.
Then the distribution $D(a,b,f,g,h)=(1,3,3,1,1)$ is not semi-greedily 
$d$-solvable, so $H$ is not semi-greedy.
Also, $H$ is not tree-solvable.
Indeed, the distribution $D(b,c,d,g,h)=(1,5,1,1,1)$ has no tree-solution
for the root $f$.
Knowing that a graph is greedy would be of great aid in solving some
pebbling distribution on the graph via computer.
Thus it is worth attempting to characterize greedy graphs.
However, this seems difficult, especially in light of the fact that
greediness is not preserved under cartesian product (nor is tree-solvability).
In fact, we will see in Section \ref{Products} that there is a graph
which is the product of two greedy graphs but which is not even
semi-greedy.

% #######################################################################
%       
%       SUBSECTION:
%       
\subsection{Diameter and Connectivity}\label{Diameter}

Another natural line of inquiry is to find necessary and sufficient conditions 
for a graph $G$ to satisfy $f(G)=n(G)$.
Although this seems very difficult, several relevant results are known.

\begin{fact}\label{Cut}
If $G$ has a cut vertex then $f(G)>n(G)$.
\end{fact}
Indeed, let $x$ be a cut vertex and $G_1,G_2$ be two components of $G-x$.
Choose $r\in V(G_1)$ and $y\in V(G_2)$.
The $r$-unsolvable distribution $D$, defined by $D(r,x,y)=(0,0,3)$ and 
$D(v)=1$ for all other $v$, witnesses this fact.

Since a graph $G$ of girth $g$ satisfies $f(G)\ge f(C_g)$, one can easily show 
the following.

\begin{fact}\label{Girth}
If $girth(G)>2\log_2 n$ then $f(G)>n(G)$.
\end{fact}

Thus it is natural to ask if there is a constant $g$ so that $girth(G)>g$ 
implies $f(G)>n=n(G)$.
More generally, is it true that for all $m$ there is a constant 
$g=g(m)$ so that $girth(G)>g$ implies $f(G)>mn$?
Both questions are completely open at the time of this writing.

Now consider upper bounds.

\begin{theorem}\cite{Snevily}\label{Diam2}
If $G$ has diameter 2 then $f(G)\le n(G)+1$.
\end{theorem}

Graphs $G$ for which $f(G)=n(G)$ are called {\it Class} 0, and otherwise
are called {\it Class} 1.
Consider the graph $H$, the union of the 6-cycle $abcdeg$ and the clique 
$ace$.
This graph has diameter two and no cut vertex, and yet $f(H)=7$ (the 
distribution $D(b,d)=(3,3)$ is $g$-unsolvable).
Clarke, Hochberg and Hurlbert \cite{CHH} characterized precisely which 
diameter two graphs are Class 1.
The characterization is based on the structure of $H$.
They developed an $O(n^5)$ algorithm which tests for Class 1 membership of 
diameter 2 graphs, and thus also serves as a test for Class 0 membership.
A key ingredient in the classification proof is showing that in any 
size $n(G)$ unsolvable distribution $D$ of a 2-connected, diameter 2, 
Class 1 graph $G$, there is no vertex containing two pebbles and there are 
exactly two vertices containing three pebbles.
Since $\max D<4$ it follows that there are at most four vertices which are 
void of pebbles, one of which must be the root $r$.
Hence, $G$ is not 4-connected, since otherwise if $D(v)=3$ then we can solve
$D$ by pebbling along one of the $vr$-paths having no void vertex.
The full characterization goes further and yields the following
result as a corollary.

\begin{theorem}\cite{CHH}\label{3Conn}
If $G$ is 3-connected and has diameter 2 then $f(G)=n(G)$.
\end{theorem}

Now, for fixed $p$ the probability that the random graph $G(n,p)$ is 
3-connected and has diameter 2 tends to 1 as $n$ tends to infinity.
Thus we have

\begin{corollary}\cite{CHH}\label{Random}
Almost all graphs are Class 0.
\end{corollary}

Since the Petersen graph $P$ is 3-connected and has diameter 2, we obtain 
a second proof that $f(P)=10$.
To extend the relationship between connectivity and diameter, Clarke et al. 
\cite{CHH} had conjectured the following, recently proved

\begin{theorem}\cite{CHKT}\label{DiamConn}
There is a function $k(d)$ such that, if $G$ is a $k(d)$-connected, 
diameter $d$ graph, then $G$ is Class 0.
\end{theorem}

The proof yields a value of $k(d)=2^{2d+3}$.
It is shown in \cite{CHH} that the function $k(d)$ must be at least 
$2^d/d$, and we believe that an optimal such $k(d)$ is less than $2^d$.
It is worth mentioning that Czygrinow et al. use Theorem \ref{DiamConn} 
to obtain the following improvement on Corollary \ref{Random}.

\begin{corollary}\cite{CHKT}\label{Threshold}
Consider the random graph $G(n,p)$ on $n$ vertices in which each edge is 
included independently with probability $p$.
Let $Q$ be the event that $G(n,p)$ is Class 0.
Then for any $d>0$ we have
\begin{description}
\item {\bf (a)}
$pn/(n\log n)^{1/d}\rar\infty$ implies that $\Pr(Q)\rar 1$ as 
$n\rar\infty$, and
\item {\bf (b)}
$pn/\log n\rar 0$ implies that $\Pr(Q)\rar 0$ as $n\rar\infty$.
\end{description}
\end{corollary}
Since being Class 0 is a monotone increasing graph property, its
threshold function exists (see \cite{BT}) and is between the two functions
above.
It would be quite interesting to find the actual value of the Class 0
threshold function.

A nice family of graphs in relation to Theorem \ref{DiamConn} is the following.
For $n\ge 2t+1$, the {\it Kneser graph}, $K(n,t)$, is the graph with
vertices ${[n] \choose t}$ and edges $\{ A,B\}$ whenever $A\cap B=\mt$.
The case $t=1$ yields the complete graph $K_n$ and the case $n=5$ and $t=2$ 
yields the Petersen graph $P$, both of which are Class 0.
When $t\ge 2$ and $n\ge 3t-1$ we have $diam(K(n,t))=2$.
Also, it is not difficult to show that $\k(K(n,t))\ge 3$ in this range,
implying that $K(n,t)$ is Class 0 by Theorem \ref{3Conn}.
Furthermore, Chen and Lih \cite{Chen} have shown that $K(n,t)$ is
connected, edge transitive, and regular of degree ${n-t\choose t}$.
A theorem of Lov\'asz \cite{Lov} states that such a graph has connectivity
equal to its degree, and thus $\k=\k(K(n,t))={n-t\choose t}$.
Therefore, using Theorem \ref{DiamConn}, it is not difficult to prove

\begin{theorem}\label{Kneser}
For any constant $c>0$, there is an integer $t_0$ such that, for $t>t_0$, 
$s\ge c(t/\log_2 t)^{1/2}$, and $n=2t+s$, we have that $K(n,t)$ is Class 0.
\end{theorem}

In the context of graph pebbling, the family of Kneser graphs is interesting 
precisely because the graphs become more sparse as $n$ decreases toward 
$2t+1$, so the diameter (as well as the girth) increases and yet the 
connectivity decreases.

\begin{question}\label{ClassKneser}
For $1\le s\ll (t/\log_2 t)^{1/2}$, is $K(2t+s,t)$ Class 0?
\end{question}
Because Pachter et al. \cite{Snevily} also proved that diameter two 
graphs have the 2-pebbling property (see below), it is interesting as well 
to ask whether $K(n,t)$ has the 2-pebbling property when $n<3t-1$ (i.e., 
when its diameter is at least 3).

% #######################################################################
%       
%       SUBSECTION:
%       
\subsection{Products and 2-Pebbling}\label{Products}

Chung proved that $f(Q^m)=2^m$.  
What she proved is, in fact, more general.
Let $\d=\la d_1,\ldots,d_m\ra$ and denote by 
$P_\d$ the graph $P_{d_1}\Box\cdots\Box P_{d_m}$.

\begin{theorem}\cite{Chung}\label{Cube}
For all nonnegative $\d=\la d_1,\ldots,d_m\ra$, we have that
$f(P_\d)$ $=2^{d_1+\cdots +d_m}$.
\end{theorem} 

The following more general conjecture has generated a great deal of
interest.

\begin{conjecture}(Graham)\label{Product}
For all $G_1$ and $G_2$, we have that $f(G_1\Box G_2)\le f(G_1)f(G_2)$.
\end{conjecture} 

There are few results which support Graham's conjecture. 
Among these, the conjecture holds for a tree times a tree \cite{Moews}, 
a cycle times a cycle (with possibly some small exceptions: it holds for 
$C_5\Box C_5$ \cite{HH}, and otherwise for $C_m\Box C_n$, provided 
$m$ and $n$ are not both from the set $\{5,7,9,11,13\}$ \cite{Snevily}), 
and a clique times a graph with the $2$-pebbling property \cite{Chung}.

A graph $G$ has the {\it 2-pebbling property} if, for any distribution $D$ of size at least $2f(G)-q(D)+1$, it is possible to move two pebbles to any
specified root $r$ after a sequence of pebbling steps.
Here, $q(D)$ is the size of the {\it support} of $D$, the number of vertices 
$v$ with $D(v)>0$.
Among the graphs known to have the 2-pebbling property are cliques, trees 
\cite{Chung}, cycles \cite{Snevily}, and diameter two graphs \cite{Snevily}.
Until recently, the only graph known not to have the 2-pebbling property 
was the {\it Lemke graph} $L$, whose vertex set is $\{a,b,c,d,w,x,y,z\}$ 
and whose edge set consists of the union of the complete bipartite graphs 
$\{a\}\times\{b,c,d\}$ and $\{b,c,d\}\times\{w,z\}$ with the path $w,x,y,z,a$.
As a witness consider that $f(L)=8$ and let 
$D(a,b,c,d,w,x,y,z)=(8,1,1,1,0,0,0,1)$.
It is impossible to move two pebbles to the root $x$.
In \cite{Foster} Foster and Snevily construct a sequence of graphs 
$L_0,L_1,L_2,\ldots$, each of which is conjectured not to have the 
2-pebbling property.
The sequence is defined starting with $L_0=L$, and $L_k$ is formed from 
$L_{k-1}$ by subdividing each of the four edges incident with $a$ exactly once.
Recently, Wang \cite{Wang} discovered a sequence of graphs 
$W_0,W_1,W_2,\ldots$ very similar to Foster and Snevily's sequence,
and proved that none of them has the 2-pebbling property.
The sequence begins with $W_0=L$, and then $W_k$ is formed from $W_{k-1}$
by subdividing each of the four edges incident with $a$ exactly once, as 
above, and then forming a clique on the 4 new vertices of the subdivisions.

The importance of the 2-pebbling property arises from its use in Chung's
proof of Theorem \ref{Cube}.
Clarke and Hurlbert \cite{Clarke,Hurl} generalize Chung's technique to cover 
a larger class of graphs than cartesian products.
Given two graphs $G_1$ and $G_2$, denote by $B(G_1,G_2)$ the set of all
bipartite graphs $F$ such that $E(F)\subseteq V(G_1)\times V(G_2)$ and
such that $F$ has no isolated vertices.
We let $\M(G_1,G_2)$ denote the set of graphs
$\{H|H=(G_1+G_2)\cup F {\rm\ for\ some\ } F\in B(G_1,G_2)\}$,
where $+$ denotes the vertex-disjoint graph union.

\begin{theorem}\cite{Clarke,Hurl}\label{GenProd}
Let $G_1$ and $G_2$ have the 2-pebbling property and $H\in \M(G_1,G_2)$.
Then $f(H)\le f(G_1)+f(G_2)$.
Furthermore, if $f(H)=f(G_1)+f(G_2)$ then $H$ has the 2-pebbling property.
\end{theorem} 

The proof follows essentially the same argument found in \cite{Chung}.
A useful corollary is that, whenever the hypotheses of Theorem \ref{GenProd}
hold, if also $f(G_i)=n(G_i)$ for each $i$ then $f(H)=n(H)$, and $H$ has the
2-pebbling property.
A pretty instance of this stronger result is a third proof that the 
Petersen graph $P$ has pebbling number 10.
Indeed, $P\in \M(C_5,C_5)$.
Moreover, we also obtain that $P$ has the 2-pebbling property because
$C_5$ does.

The power of Theorem \ref{GenProd} is shown in that it is used easily 
to prove its own generalization.
Denote by $\M(G_1,\ldots,G_t)$ the set of all graphs $H$ such that
$H[V_i\cup V_j]\in \M(G_i,G_j)$ for all $i\ne j$, where $V_i=V(G_i)$.
Here, $H[X]$ denotes the subgraph of $H$ induced by the vertex set $X$.

\begin{theorem}\cite{Clarke,Hurl}\label{GenProd2}
Let $G_i$ have the 2-pebbling property for $1\le i\le t$ and let
$H\in \M(G_1,\ldots,G_t)$.
Then $f(H)\le \sum\limits_{i=1}^tf(G_i)$.
Furthermore, if $f(H)=\sum\limits_{i=1}^tf(G_i)$ then $H$ has the 2-pebbling 
property.
\end{theorem} 

The analogous corollary is that, whenever the hypotheses of Theorem 
\ref{GenProd2} hold, if also $f(G_i)=n(G_i)$ for each $i$ then $f(H)=n(H)$, 
and $H$ has the 2-pebbling property.
This corollary yields a simple proof of another result of Chung, and
verifies another instance of Graham's conjecture.
Notice that $G\Box K_m\in \M(G,\ldots,G)$.

\begin{theorem}\cite{Chung}\label{CliqueProd}
If $G$ has the 2-pebbling property then $f(G\Box K_m)\le mf(G)$.
\end{theorem} 

A similar result was proved by Moews:

\begin{theorem}\cite{Moews}\label{TreeProd}
If $G$ has the 2-pebbling property and $T$ is a tree then 
$f(G\Box T)\le f(G)f(T)$.
\end{theorem} 
Since trees have the 2-pebbling property, Theorem \ref{TreeProd} shows
that the cartesian product of two trees also satisfies Graham's conjecture.

Finally, we remark that the pebbling number of a product can sometimes fall
well inside the range $n(G_1\Box G_2)<f(G_1\Box G_2)<f(G_1)f(G_2)$.
Consider the graph $H=P_3\Box S_4$.
($S_n$ is the {\it star} with $n$ vertices, also denoted $K_{1,n-1}$).
It is easy to see that $f(P_3)=4$ and $f(S_4)=5$.
Although each pebbling number is only one more than the corresponding
number of vertices, $f(H)=18$ is far greater than $n(H)=12$; it would be 
interesting to investigate how much bigger this gap can be in general.
Also, notice that $f(H)<f(P_3)f(S_4)$, a strict inequality.
More importantly, as observed by Moews \cite{M}, $H$ is not semi-greedy.
Indeed, think of $H$ as three pages of a book, let $r$ be the corner vertex
of one of the pages, $x$ the farthest corner vertex of a second page, 
$u$, $v$ and $w$ the three vertices of the third page, and let
$D(u,v,w,x)=(1,1,1,15)$.
This shows that even semi-greediness is not always enjoyed by the product 
of two greedy graphs.
This is especially disappointing since it shoots down a promising method 
of attack on Graham's conjecture.

% #######################################################################
% #######################################################################
%       
%       NEW SECTION:
%       
\section{Number Theory}\label{Number}
The concept of pebbling in graphs arose from an attempt by Lagarias and
Saks \cite{Saks} to give an alternative proof of a theorem of Kleitman 
and Lemke.
An elementary result in number theory which follows from the pigeonhole
principle is 

\begin{fact}\label{Pigeon}
For any set $N=\{n_1,\ldots,n_q\}$ of $q$ natural numbers,
there is a nonempty index set $I\subset\{1,\ldots,q\}$ such that
$q\ \Big|\sum\limits_{i\in I}n_i$.
\end{fact}

Erd\H os and Lemke conjectured in 1987 that the extra condition
$\sum\limits_{i\in I}n_i\le {\rm lcm}(q,n_1,\ldots,n_q)$ could also be
guaranteed.
In 1989 Lemke and Kleitman proved

\begin{theorem}\cite{KL}\label{KleitLem}
For any any set $N=\{n_1,\ldots,n_q\}$ of $q$ natural numbers,
there is a nonempty index set $I\subset\{1,\ldots,q\}$ such that
$q\ \Big|\sum\limits_{i\in I}n_i$ and 
$\sum\limits_{i\in I}{\rm\ gcd}(q,n_i)\le q$.
\end{theorem} 
This proves the Erd\H os-Lemke conjecture because of the string of
inequalities 
$$\sum\limits_{i\in I}n_i
={1\over q}\sum\limits_{i\in I}qn_i
={1\over q}\sum\limits_{i\in I}{\rm\ lcm}(q,n_i){\rm\ gcd}(q,n_i)$$
$$\le {1\over q}{\rm\ lcm}(q,n_1,\ldots,n_q)
\sum\limits_{i\in I}{\rm\ gcd}(q,n_i)
\le {\rm\ lcm}(q,n_1,\ldots,n_q).$$

The argument used by Kleitman and Lemke had many cases and did not seem
to be the most natural proof.
It was the intention of Lagarias and Saks to introduce graph pebbling as
a more intuitive vehicle for proving the theorem.
If the formula for the general pebbling number of a cartesian product of paths 
is as was believed, then the number-theoretic result would follow easily.
It was Chung \cite{Chung} who finally pinned down such a formula.
Recently, pebbling has been used to extend the result further.
Denley \cite{Denley} proved that if each $n_i|q$ (with $n_i\le n_{i+1}$) 
and $\sum\limits_{p{\rm\ prime,\ }p|q}1/p\le 1$, then there is a nonempty 
$I$ such that $q=\sum\limits_{i\in I}n_i$ and $n_i|n_j$ for all $i<j$.

Kleitman and Lemke went on to make more general conjectures on groups.
First, let $G$ be a finite group of order $q$ with identity $e$, and let
$|g|$ denote the order of the element $g$ in $G$.
Then for any multisubset $N=\{g_1,\ldots,g_q\}$ of $G$ there is a
nonempty $I$ such that $\prod\limits_{i\in I}\ g_i=e$ and 
$\sum\limits_{i\in I}1/|g_i|\le 1$.
Their prior theorem is merely the case $G=\Z_q$, and they verified this
conjecture for $G=\Z_p^n$, for dihedral $G$, and also for
all $q\le 15$.
Second, let $H$ be a subgroup of a group $G$ with $|G/H|=q$ and let
$N=\{g_1,\ldots,g_q\}$ be any multisubset of $G$.
Then there is a nonempty $I$ such that $\prod\limits_{i\in I}\ g_i\in H$ and
$\sum\limits_{i\in I}1/|g_i|\le 1\Big/|\prod\limits_{i=1}^q\ g_i|$.
The first conjecture is the case $H=\{e\}$ here, and this second conjecture
they verified for all $|G|\le 11$.
It would be interesting to see what pebbling could say about these
two conjectures.

In order to describe Chung's proof (see \cite{Chung,CHH}) of Theorem 
\ref{KleitLem} we need to define a more general pebbling operation on a
product of paths.

A {\it $p$-pebbling step} in $G$ consists of removing $p$ pebbles from a
vertex $u$, and placing one pebble on a neighbor $v$ of $u$.  
The definitions for $r$-solvability, and so on, carry over to $p$-pebbling.
Recall the definition of the graph $P_\d$ from Section \ref{Products}.
Each vertex $v\in V(P_\d)$ can be represented by a vector 
$\vv=\la v_1,\ldots,v_m\ra$, with $0\le v_i\le d_i$ for each $i\le m$.
Let $\e_i=\la 0,\ldots,1,\ldots,0\ra$ be the $i^{\rm th}$ standard 
basis vector and $\0=\la 0,\ldots,0\ra$.
Then two vertices $u,v$ are adjacent in $P_\d$ if and 
only if $\u-\vv=\pm\e_i$ for some integer $1\le i\le m$.
If $\p=(p_1,\ldots,p_m)$, then we define $\p$-{\it pebbling} in 
$P_\d$ to be such that each pebbling step from $\u$ to $\vv$ is a 
$p_i$-pebbling step whenever $\u-\vv=\pm\e_i$.
The $\p$-pebbling number of $P_\d$ is denoted by $f_\p(P_\d)$.

For integers $p_i,d_i\ge 1$, $1\le i\le m$, we use $\p^\d$
as shorthand for the product $p^{d_1}_1\cdots p^{d_m}_m$.
Chungs's proof uses the following result.

\begin{theorem}\cite{Chung}\label{GenPathProd}
Every distribution of size at least $\p^\d$ is $\0$-solvable via 
greedy $\p$-pebbling.
\end{theorem}

Actually one can prove more.

\begin{theorem}\cite{CHH}\label{GenPebNum}
The $\p$-pebbling number $f_\p(P_\d)=\p^\d$.  
Moreover, $P_\d$ is greedy.
\end{theorem}

In order to prove Theorem \ref{KleitLem} from Theorem \ref{GenPathProd}
one first defines a pebbling distribution $D$ in $P_\d$ which depends 
on the set of integers $\{x_1,\ldots,x_d\}$.
Here, $|D|=d=\p^\d =\prod\limits_{i=1}^{m}p_i^{d_i}$, the prime 
factorization of $d$, where $\d=\la d_1,\ldots,d_m\ra$.
In what follows, each pebble will be named by a set, and $\c(B)$
will denote the vertex (coordinates) on which the pebble $B$ sits.
We let $x_j$ correspond to the pebble $A_j=\{x_j\}$ , which we place on the
vertex $\c(A_j)=\la c_1,\ldots,c_m\ra$ of $P_\d$, where 
$d/{\rm gcd}(x_j,d)=\p^\c$.
For each vertex $\u=\la u_1,\ldots,u_m\ra$ define the set
$X(\u)=\{A|\c(A)=\u\}$ to denote those pebbles currently sitting on 
$\u$, and let $\u^{(i)}=\la u_1,\ldots,u_{i}-1,\ldots,u_m\ra$.

For a set $B$ we make the following recursive definitions.
The {\it value} of $B$ is defined as $val(B)=\sum\limits_{A\in B}val(A)$,
with $val(\{A_j\})=x_j$.
The function $GCD$ is defined as $GCD(B)=\sum\limits_{A\in B}GCD(A)$,
where $GCD(\{A_j\})={\rm gcd}(x_j,d)$.
Finally, $Set(B)=\bigcup\limits_{A \in B}Set(A)$, where $Set(A_j)=A_j$.

We say that $B$ is {\it well placed} at $\c(B)=\la c_1,\ldots,c_m\ra$ when
\begin{equation}\label{WP1}
\p^{\d-\c(B)}\vert val(B)
\end{equation}
and
\begin{equation}\label{WP2}
GCD(B)\le \p^{\d-\c(B)}\ .
\end{equation}

It is important to maintain a numerical interpretation of $\p$-pebbling 
so that moving a pebble to $\0$ corresponds to finding a set $J$ which 
--- playing the role of $I$ --- satisfies the conclusion of 
Theorem \ref{KleitLem}.
For this reason we introduce the following operation, which corresponds to
a greedy $p_i$-pebbling step in which a numerical condition must hold in
order to move a pebble.
It is shown that this condition holds originally for $D$ (Lemma \ref{StartWP}) 
and is maintained throughout (Lemma \ref{KeepWP}).

\begin{quote}
{\bf Numerical Pebbling Operation.}
If $W$ is a set of $p_i$ pebbles such that every pebble $A\in W$ sits on 
the vertex $\c(A)=\u$, and there is some $B\subseteq W$ such that 
$p_i^{b_i}|val(B)$, where $b_i=d_i-c_i+1$, then replace $X(\c)$ by 
$X(\c)\setminus W$, and replace $X(\c^{(i)})$ by $X(\c^{(i)})\cup B$.
\end{quote}

\noindent
\begin{lemma}\label{StartWP}
$A_j$ is well placed for $1\le j\le d$.
\end{lemma}

\begin{lemma}\label{KeepWP}
Suppose $B \subseteq X(\u), |B|\le p_i$, and $p_i^{b_i}|val(B)$ for 
$b_i=d_i-u_i+1$.
Suppose further that for every $A\in B$, $A$ is well placed at $\u$.
Then $B$ is well placed at $\u^{(i)}$.
\end{lemma}

\begin{lemma}\label{equiv}
Suppose $|X(\u)|\ge p_i$, and for all $A\in X(\u)$, $A$ is well placed 
at $\u$.
Then there exists some $B\subseteq X(\u)$ such that $|B|\le p_i$ and 
$p_i^{b_i}|val(B)$ where $b_i=d_i-u_i+1$.
\end{lemma}

By Lemma \ref{StartWP} the pebbles corresponding to each of the numbers 
are initially well placed.
Lemma \ref{KeepWP} guarantees that applying the Numerical Pebbling Operation 
maintains the well placement of the pebbles.
Lemma~\ref{equiv} establishes that every graphical pebbling operation 
can be converted to a numerical pebbling operation.
Then by Theorem \ref{GenPathProd} we can repeatedly apply the numerical 
pebbling operation to move a pebble to $\0$.
This pebble $B$ is then well placed at $\0$.
Thus, for $J=\{ j|x_j\in Set(B)\}$, we have
$d=\p^\d\ |\ val(B)=\sum\limits_{j\in J}x_j$ by (\ref{WP1}), and 
$\sum\limits_{j\in J}gcd(x_j,d)=\sum\limits_{j\in J}x_j=GCD(B)\le
\p^\d=d$ by (\ref{WP2}).
This proves Theorem \ref{KleitLem}.
Interestingly, it is Fact \ref{Pigeon} which is used to prove
Lemma \ref{equiv}.

A natural generalization of Graham's Conjecture \ref{Product} is
the following 

\begin{conjecture}\cite{CHH}\label{pProduct}
For all $G_1, G_2, p_1, p_2$, we have that
$f_{(p_1,p_2)}(G_1\Box G_2)\le$ $f_{p_1}(G_1)f_{p_2}(G_2)$.
\end{conjecture} 

We leave it for the reader to ponder whether this is the ``right''
generalization, say, for driving an inductive proof of Conjecture 
\ref{Product}.

% #######################################################################
% #######################################################################
%       
%       NEW SECTION:
%       
\section{Thresholds\protect\footnotemark[2]}\label{Thresholds}
\addtocounter{footnote}{2}
\footnotetext{Much of this section is adapted from \protect\cite{CEHK}.}
In this section, we introduce a probabilistic pebbling model,
where the pebbling distribution is selected uniformly at random from
the set of all distributions with a prescribed number of pebbles.
(Here the distribution of pebbles to vertices is like the distribution
of unlabeled balls to labeled urns.)
We define and study thresholds for the number $t$ of pebbles so that if
$t$ is essentially larger than the threshold, then any
distribution is almost surely solvable, and if $t$ is essentially
smaller than the threshold, then any distribution is almost surely unsolvable.
Of course, the definition mimics the important threshold concept in random
graph theory.
Unlike the situation in random graphs, however, it does not seem obvious
that even ``natural'' families of graphs have pebbling thresholds.
One candidate for such a family is the sequence of paths.
We should emphasize also that, unlike in random graph theory, even the
most basic random variables considered here are functions of dependent
random variables, and the dependence is not ``sparse''.
This substantially limits the set of tools available for analyzing these
random variables.

We now recall some basic asymptotic notation.
For two functions $f$ and $g$, we write $f\ll g$ (equivalently $g\gg f$)
when the ratio $f(n)/g(n)$ approaches $0$ as $n$ tends to infinity.
We use $o(g)$ and $\w(f)$, respectively, to denote the sets
$\{f\ |\ f\ll g\}$ and $\{g\ |\ f\ll g\}$,
so that $f\in o(g)$ if and only if $g\in\w(f)$.
In addition, we write $f\in O(g)$ (equivalently $g\in\W(f)$) when there
are positive constants $c$ and $k$ such that $f(n)/g(n)<c$ for all $n>k$,
and we write $\Th(g)$ for $O(g)\cap\W(g)$.
We also use the shorthand notation $\Th(f)\le \Th(g)$ to mean that
$f\pr\in O(g\pr)$ for every $f\pr\in \Th(f)$ and $g\pr\in \Th(g)$.
To avoid cluttering the paper with floor and ceiling symbols, we adopt the
convention that large constants (such as $1/\ep$ when $\ep$ is small)
are integers.

We are almost ready to define formally our notion of a pebbling
threshold function.
Let $D_n:[n]\rar \N$ denote a distribution of pebbles on $n$ vertices.
For a particular function $t = t(n)$, we consider the probability space
$\Omega_{n,t}$ of all distributions $D_n$ of {\it size} $t$, i.e. with
$t=\sum_{i\in[n]}D_n(i)$ pebbles.
Given a graph sequence $\G=(G_1,\ldots,G_n,\ldots)$, denote by $P_\G(n,t)$
the probability that an element of $\Omega_{n,t}$ chosen uniformly at
random is $G_n$-solvable.
We call a function $g$ a {\it threshold} for $\G$, and write $g\in th(\G)$,
if the following two statements hold as $n\rar\inf$:
({\it i}) $P_\G(n,t)\rar 1$ whenever $t\gg g$, and
({\it ii}) $P_\G(n,t)\rar 0$ whenever $t\ll g$.

We shall consider the following families of graphs.
\begin{itemize}
\item $\K=(K_1,\ldots,K_n,\ldots)$: $K_n$ is the complete graph on $n$ vertices.
\item $\P=(P_1,\ldots,P_n,\ldots)$: $P_n$ is the path on $n$ vertices.
\item $\C=(C_1,\ldots,C_n,\ldots)$: $C_n$ is the cycle on $n$ vertices.
\item $\S=(S_1,\ldots,S_n,\ldots)$: $S_n$ is the star on $n$ vertices.
\item $\Wh=(W_1,\ldots,W_n,\ldots)$: $W_n$ is the wheel on $n$ vertices.
\item $\Q=(Q^1,\ldots,Q^m,\ldots)$: $Q^m$ is the $m$-dimensional cube
on $n=2^m$ vertices.
\end{itemize}
In addition, we will denote generic families of graphs by
$\G=(G_1,\ldots,G_n,\ldots)$ or $\cH=(H_1,\ldots,H_n,\ldots)$.

Although the existence of $th(\G)$ has yet to be established and may be 
impossible for some graph sequences $\G$, the family of cliques
is not one of them.
(If $th(\G)\ll th(\cH)$ then the alternating graph sequence
$(G_1,H_2,G_3,H_4,\ldots)$ may seem to have no threshold function
--- c.f. Theorems \ref{CliqueThresh} and \ref{Thresh}(f) ---
its threshold function is merely the alternation of the corresponding two
threshold functions.
There is no demand on the continuity of such functions.)

\begin{theorem}\cite{Clarke}\label{CliqueThresh}
The threshold for cliques is $th(\K)=\Th(\sn)$.
\end{theorem}

This result is merely a reformulation of the so-called ``Birthday problem''
\cite{Feller} in which one finds the probability that 2 of $t$ people 
share the same birthday, assuming $n$ days in a year.

Among the results of Czygrinow et al. are the following.

\begin{theorem}\cite{CEHK}\label{Thresh}
\begin{description}
\item {\bf (a)}
For all $\G$, if $th(\G)$ exists then, for every $\ep>0$, we have
$th(\G)\subseteq \Omega(\sn)\cap o(n^{1+\ep})$.
\item {\bf (b)}
Let $d(n)=diameter(G_n)$ and suppose that $th(\G)$ exists.
If $d(n)\le d$ for all $n$, then $th(\G)\subseteq O(n)$.
\item {\bf (c)}
Let $d(n)=diameter(G_n)$, $k(n)=connectivity(G_n)$, and suppose that 
$th(\G)$ exists.
If $k(n)\ge 2^{2d(n)+3}$ for all $n$, then $th(\G)\subseteq O(n)$.
\item {\bf (d)}
If $th(\Q)$ exists then $th(\Q)\subseteq \Omega(\sn)\cap O(n)$.
\item {\bf (e)}
If $th(\C)$ exists then, for every $\ep>0$, we have 
$th(\C)\subseteq \W(n)\cap o(n^{1+\ep})$.
\item {\bf (f)}
If $th(\P)$ exists then, for every $\ep>0$, we have 
$th(\P)\subseteq \W(n)\cap o(n^{1+\ep})$.
\item {\bf (g)}
$th(\S)=\Th(\sn)$.
\item {\bf (h)}
$th(\Wh)=\Th(\sn)$.
\end{description}
\end{theorem}
Two statements within Theorem \ref{Thresh} are somewhat surprising,
namely, that no tighter bounds for cubes or paths are known than
what is found in (d) and (f).
Also it is interesting to note that, in light of the results for paths (f)
and stars (g), it is conceivable that the set of thresholds for all possible 
sequences of trees may span the entire range of functions from $n^{1/2}$ 
to $n$ (or $n^{1+\ep}$, as the case may be).
As witnessed by Theorem \ref{Thresh} (b,c), diameter seems to be a critical 
parameter.
In order to make some of these remarks more precise, we offer the following

\begin{conjecture}\label{Range}
For every $t_1$ and $t_2$ such that $t_1\in\Omega(\sn)$, $t_2\in O(n)$, 
and $t_1\ll t_2$,  there is a graph sequence $\G=\{G_1,\ldots,G_n,\ldots\}$ 
such that $th(\G)\in\Omega(t_1)\cap O(t_2)$.
Moreover, there is such a sequence in which each $G_n$ is a tree.
\end{conjecture}

To date the theory of graph pebbling omits two important assertions, namely
an existence theorem for threshold functions of certain families 
of graph sequences and a monotonicity theorem for threshold functions.

\begin{conjecture}\label{ExistThresh}
Every graph sequence $\G=\{G_1,\ldots,G_n,\ldots\}$ has a 
threshold $th(\G)$.
\end{conjecture}

\begin{conjecture}\label{MonoThresh}
If $f(G_n)\le f(H_n)$ for all $n$ and both $th(\G)$ and $th(\cH)$ exist,
then $th(\G)\le th(\cH)$.
\end{conjecture}

For a positive integer $t$ and a graph $G$ denote by $p(G,t)$ the probability 
that a randomly chosen distribution $D$ of size $t$ on $G$ solves $G$.
If both $th(\G)$ and $th(\cH)$ exist, then Conjecture \ref{MonoThresh}
would follow from the statement that, if $f(G_n)\le f(H_n)$ then for
all $t$ we have $p(G_n,t)\ge p(H_n,t)$.
Unfortunately, although seemingly intuitive, this implication is false.
Using the Class 0 pebbling characterization theorem of \cite{CHH}, we
discovered in \cite{CEHK} a family of pairs of graphs $(G_n,H_n)$, one pair
for each $n=3k+4$, for which the implication fails.
Figure~1 suggests the structure of these graphs; see \cite{CEHK} for details.

\begin{figure}
\begin{multicols}{2}
\centerline{\hbox{\psfig{figure=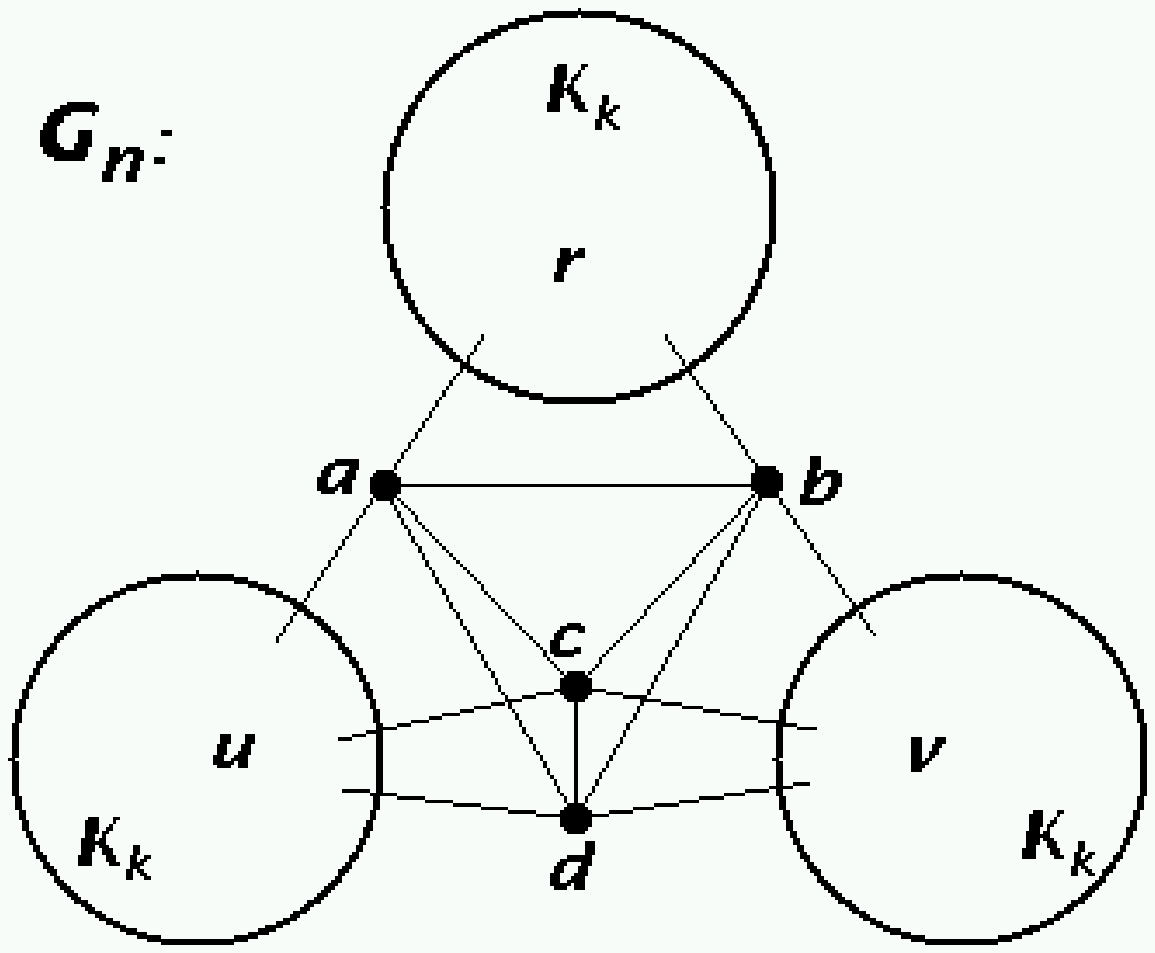,height=2.2in}}}

\centerline{\hbox{\psfig{figure=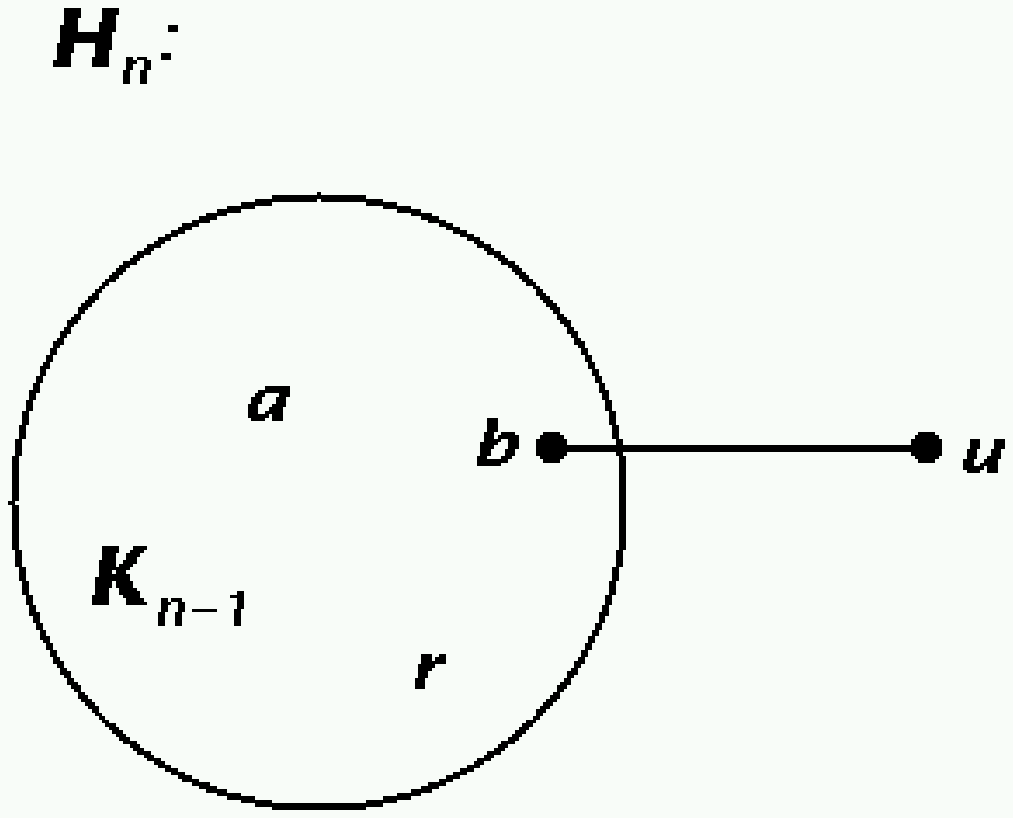,height=2.2in}}}
\end{multicols}
\caption{A counterexample to $f(G_n)\le f(H_n)\Rightarrow 
p(G_n,t)\ge p(H_n,t)$.}
\end{figure}

% #######################################################################
% #######################################################################
%       
%       NEW SECTION:
%       
\section{Set Theory}\label{Sets}
In order to describe the connection between extremal set theory and 
threshold functions we need to introduce some notation.

Let $[n]=\{1,2,\ldots,n\}$,
$Set\{n\}$ be the partially ordered set ({\it poset}) of subsets of $[n]$,
$MSet\{n\}$ be the poset of submultisets of $[n]$, and
$BMSet\{n,b\}$ be the poset of $b$-bounded submultisets of $[n]$
(that is, no element appears more than $b$ times).
Ordered by the relation of inclusion, we call $Set\{n\}$ the
{\it subset lattice} (also called the Boolean algebra) and $BMSet\{n,b\}$ 
the {\it multiset lattice} (the product of $n$ chains $C_{b+1}$ on 
$b+1$ elements each; isomorphic to the lattice of divisors of an integer
which is the $b^{\rm th}$ power of a product of $n$ distinct primes).

Let $Set[n,w]$ denote the family of weight-$w$ subsets of $[n]$
(i.e., $Set[n,w]={[n]\choose w}$),
$MSet[n,w]$ denote the family of weight-$w$ submultisets of $[n]$, and
$BMSet[n,w,b]$ denote the family of $b$-bounded weight-$w$ submultisets 
of $[n]$.
Also let 
\begin{center}
\begin{tabular}{rcl}
$bin[n,w]$&$=$&$\Big|Set[n,w]\Big|={n\choose w}$,\\
\ \\
$mul[n,w]$&$=$&$\Big|MSet[n,w]\Big|={n+w-1\choose w}$, and\\
\ \\
$bmul[n,w,b]$&$=$&$\Big|BMSet[n,w,b]\Big|$\\
\ \\
	&$=$&$\sum\limits_{i=0}^{\lf w/(b+1)\rf} 
		(-1)^i\ bin[n,i]\ mul[n,w-i(b+1)]$.
\end{tabular}
\end{center}

Next, let $\N =\{1,2,\ldots\}$,
$S[w]$ be the set of weight-$w$ subsets of $\N$,
$MS[w]$ be the set of weight-$w$ submultisets of $\N$, and
$BMS[w,b]$ be the set of $b$-bounded weight-$w$ submultisets of $\N$.
Since $BMS[w,b]$ is the most general setting ($S[w]=BMS[w,1]$ and
$MS[w]=BMS[w,w]$), we will define the well-known {\it Colex} order $L$
on $BMS[w,b]$.

Let $M\in BMS[w,b]$ with multiplicities $(m_1,m_2,\ldots,m_l)$ 
for some $l$ so that $m_j=0$ for all $j>l$.
Define the function $C:BMS[w,b]\rar\N$ by
$C(M)=\sum\limits_{i=1}^l m_i(b+1)^{i-1}$.
Then for $A,B\in BMS[w,b]$ we have $L(A,B)$ if and only if $C(A)<C(B)$.
It is easy to see that, for any $f$ and $n$, the first $f$ multisets of
$BMSet[n,w,b]$ are precisely the first $f$ members of $BMS[w,b]$.

To define {\it shadows}, again let $M$ be an element of $BMS[w,b]$ with 
multiplicities $(m_1,m_2,\ldots,m_l)$ for some $l$ so that $m_j=0$ for all 
$j>l$.
Now denote by $Shad[M]$ the set of all $A\in BMS[w-1,b]$ with multiplicities
$(a_1,a_2,\ldots,a_l)$ so that $a_i\le m_i$ for all $i\le l$.
When $\F$ is a family of multisets, each of which is in $BMS[w,b]$,
denote by $Shad[\F]$ the family of all $A\in Shad[M]$ for some $M\in\F$,
and let $shad[\F]=\Big|Shad[\F]\Big|$.
For the vector $\vv =\la v_1,v_2,\ldots,v_s\ra$, let $r$ be such that
$v_i=0$ for all $i<r$, and define $w_j=\sum_{i=1}^j v_i$ and $w=w_s$.
Then define 
$col[\vv,w,b]=\sum_{j=r+1}^s\sum_{i=1}^{v_j}bmul[j-1,w_j+1-i,b]$
and let $Col[\vv ,w,b]$ be the first $col[\vv ,w,b]$ multisets in the 
Colex order on $BMS[w,b]$.
Finally, define $Shad[\vv ,w,b]=Shad[Col[\vv ,w,b]]=Col[\vv ,w-1,b]$
and let $shad[\vv ,w,b]=\Big|Shad[\vv ,w,b]\Big|=col[\vv ,w-1,b]$.
It is not difficult to show that for any natural numbers $f$, $w$ and $b$,
there exist natural numbers $s,v_1,v_2,\ldots,v_s$ such that $f=col[\vv ,w,b]$.
In 1969 Clements and Lindstr\" om proved the following

\begin{theorem}\cite{CL}\label{ClemLind}
Let $\vv$, $w$ and $b$ be given, with $\F\in BMS[w,b]$ and 
$|\F |=col[\vv ,w,b]$.
Then $shad[\F]\ge shad[\vv ,w,b]=col[\vv ,w-1,b]$.
\end{theorem}
The important special cases when $b=1$ and $b\ge w$ had been proven earlier
by Kruskal \cite{Kr} and Katona \cite{Kat} and by Macaulay
\cite{Mac}, respectively.
The Kruskal-Katona theorem plays a crucial role in many combinatorial contexts,
including the dimension theory of partially ordered sets (\cite{HKT}).
It is of fundamental importance in the probability estimates below.
\vspace{0.2 in}

\noindent
{\bf Note.}
Theorem \ref{ClemLind} is actually one instance of the 
Clements-Lindstr\" om theorem.
It is easy to generalize the setting of multisets uniformly bounded by $b$
to coordinate-wise bounded by $\bb =\la b_1,b_2,\ldots,b_n\ra$, as follows.
We say $m\in BMSet[n,w,\bb]$ if $M$ has multiplicities $M=(m_1,\ldots,m_n)$ 
with $0\le m_i\le b_i$ for all $i$, and we define 
$C(M)=\sum_{i=1}^nm_i\prod_{j=1}^{i-1}(b_i+1)$.
With shadows and functions such as $col[\vv,w,\bb]$ defined in the obvious 
ways, Clements and Lindstr\" om proved the analogous generalization of 
Theorem \ref{ClemLind}.
Nevertheless, we will stick to the uniform case here.
\vspace{0.2 in}

Since $bmul[n,w,b]$ is a polynomial in the variable $n$, it is well-defined
even when $n$ takes on some real value $x$.
Of course, for any natural numbers $f$, $w$ and $b$, there is a real number
$x$ for which $f=bmul[x,w,b]$.
In 1979 Lov\'asz proved the following version of the Kruskal-Katona theorem.

\begin{theorem}\cite{Lov}\label{Lovasz}
Let $w$ be given with $\F\in BMS[w,1]=S[w]$ and 
$|\F |=bmul[x,w,1]$ $=bin[x,w]$.
Then $shad[\F]\ge bmul[x,w-1,1]=bin[x,w-1]$.
\end{theorem}

For a general ranked poset $\P$ and family $\F$ of elements of $\P$, 
let $\F_w$ be those elements of $\F$ of rank $w$.
Then define the probability $p(\F_w)=|\F_w|\Big/|P_w|$.
A family $\F$ is {\it monotone decreasing} if $A\subset B\in\F$ implies
$A\in\F$.
Also, $\F$ is an {\it antichain} if no pair of elements of $\F$ are 
related in $\P$.
The poset $\P$ is {\it LYM} (has the {\it LYM property}) if, for any
antichain $\F$ in $\P$, $\sum_wp(\F_w)\le 1$.
Also, $\P$ is {\it normal} (has the {\it normalized matching property}) if,
for any monotone decreasing family $\F$, we have $p(\F_u)\ge p(\F_w)$
whenever $0<u<w$.
Finally, we define $\P$ to be {\it supernormal} if $p(\F_u)^w\ge p(\F_w)^u$
whenever $0<u<w$.

It is known that $\P$ is LYM if and only if $\P$ is normal, and that 
$Set\{n\}$ is LYM, and hence normal, for all $n$ \cite{Bol,Lub,Mesh,Yam}.
The product theorems of Canfield \cite{Can} and Harper \cite{Harper} show that 
$BMSet\{n,b\}$ is LYM, and hence normal \cite{BezEng,Hsieh}, for all $b\le n$.
An important consequence of Theorem \ref{Lovasz} is that $Set\{n\}$ is
supernormal.
It is precisely this inequality which allows one to prove the following

\begin{theorem}\cite{BT}\label{BollThom}
If $\F$ is a monotone decreasing family of subsets of $[n]$,
then there is a threshold $th(\F)$ for $\F$;
that is, $p(\F_w)\rar 1$ when $w\ll th(\F)$ and
$p(\F_w)\rar 0$ when $w\gg th(\F)$.
\end{theorem}
Of course, the analogous theorem holds for monotone increasing families.
Most notably, as a corollary one obtains the existence for threshold
functions of monotone properties of graphs (such as for connectedness,
hamiltonicity, or subgraph containment, but not for induced subgraph
containment).

In our case we would like to mimic these results for multisubsets in order
to prove existence for pebbling thresholds of graph sequences.
Unfortunately, not all of the results generalize.
We have discovered computationally that the Lov\'asz-type version (Theorem 
\ref{Lovasz}) of both the Macauley and Clements-Lindstr\"om theorems 
may hold in general.
More precisely, we make the following

\begin{conjecture}\label{GenLov}
Let $w$ and $b>1$ be given with $\F\in BMS[w,b]]$ and let $x$ be
defined by $|\F |=bmul[x,w,b]$.
Then $shad[\F]\ge bmul[x,w-1,b]$.
\end{conjecture}

Although the truth of Conjecture \ref{GenLov} would surely have applications
elsewhere, our attempt to generalize Theorem \ref{Lovasz} is motivated by our
desire to prove that $BMSet\{n,b\}$ is supernormal for all $n$ and $b$.
However, 

\begin{theorem}\label{NotSuper}
For all $n$ and $b$, $BMSet\{n,b\}$ is not supernormal.
\end{theorem}

\noindent
{\it Proof}.
If $BMSet\{n,b\}$ is supernormal then it should be the case that, for all
$\F$, $p(\F)^{b-1}<p(Shad[\F])^b$.
However, let $\vv$ be the vector having $r-1$ zeroes followed by a single $b$,
so that $r=s<n$, and consider the family $\F=Col[\vv,b,b]$.
Then $|\F|=col[\vv,b,b]={s+b-1\choose b}$ and 
$shad[\F]=col[\vv,b-1,b]={s+b-2\choose b-1}$.
Also, $bmul(n,b,b)={n+b-1\choose b}$ and $bmul(n,b-1,b)={n+b-2\choose b-1}$,
so that $p(\F)={s+b-1\choose b}\Big/{n+b-1\choose b}$ and
$p(Shad[\F])={s+b-2\choose b}\Big/{n+b-2\choose b}$.
Thus 

\medskip
$p(\F)^{b-1}-p(Shad[\F])^b$
$$\qquad =\ \left[{(s+b-2)\cdots (s)\over (n+b-2)\cdots (n)}\right]^{b-1}
	\left[\left({s+b-1\over n+b-1}\right)^{b-1}
		-{(s+b-2)\cdots (s)\over (n+b-2)\cdots (n)}\right]\ ,$$
which is positive because $x/y>(x-i)/(y-i)$ whenever $i<x<y$.
\hfill$\Box$

Our purpose for investigating the supernormality of $BMSet\{n,b\}$ derives
from our attempt to generalize Theorem \ref{BollThom}, which would yield as
a corollary the existence of pebbling threshold functions for arbitrary graph 
sequences.
Conceivably, such a generalization may be true even in the absence of 
supernormality.
For example, empirical evidence seems to indicate that, in the 
Clements-Lindstr\"om setting, for fixed $u$ and $w$, 
$\Big|p(\F_u)^w-p(\F_w)^u\Big|\rar 0$ as $n\rar\infty$.
Thus one could closely approximate $p(\F_u)^w$ by $p(\F_w)^u$.
Unfortunately, in order to prove Conjecture \ref{GenLov}, one needs such 
approximations for all $u,w$ and $n$.

% #######################################################################
% #######################################################################
%       
%       NEW SECTION:
%       
\section{Acknowledgments}\label{Ack}
The author would like to thank Mark Kayll for his generous reading of an
early draft of this paper, and for his invaluable comments.

% #######################################################################
% #######################################################################
%       
%       BIBLIOGRAPHY:
%       
\bibliographystyle{plain}

\begin{thebibliography}{999}
\bibitem{BezEng} 
S. Bezrukov and K. Engel,
{\it Properties of graded posets preserved by some operations},
The mathematics of Paul Erd\H os, II, 
Algorithms Combin. {\bf 14} (1997), 79--85.
\bibitem{Bol} 
B. Bollob\'as,
{\it Combinatorics},
Cambridge Univ. Press, Cambridge, New York (1986).
\bibitem{BT} 
B. Bollob\'as and A. Thomason,
{\it Threshold functions},
Combinatorica {\bf 7} (1987), 35--38.
\bibitem{Can} 
E.R. Canfield,
{\it A Sperner property preserved by products},
Linear and Multilinear Alg. {\bf 9} (1980), 151--157.
\bibitem{Chen} 
B.L. Chen and K.W. Lih,
{\it Hamiltonian uniform subset graphs},
J. Combin. Theory (Ser. B) {\bf 42} (1987), 257--263.
\bibitem{Chung} 
F.R.K. Chung,
{\it Pebbling in hypercubes},
SIAM J. Disc. Math. {\bf 2} (1989), 467--472.
\bibitem{Clarke} 
T.A. Clarke,
{\it Pebbling on graphs},
Master's Thesis, Arizona St. Univ. (1996).
\bibitem{CHH} 
T.A. Clarke, R.A. Hochberg and G.H. Hurlbert,
{\it Pebbling in diameter two graphs and products of paths},
J. Graph Th. {\bf 25} (1997), 119--128.
\bibitem{CL} 
G.F. Clements and B. Lindstr\" om,
{\it A generalization of a combinatorial theorem of Macaulay},
J. Combin. Theory {\bf 7} (1969), 230--238.
\bibitem{CHKT} 
A. Czygrinow, G.H. Hurlbert, H.A. Kierstead, and W.T. Trotter,
{\it A note on graph pebbling},
preprint (1999).
\bibitem{CEHK} 
A. Czygrinow, N. Eaton, G.H. Hurlbert, and P.M. Kayll,
{\it On pebbling threshold functions for graph sequences},
preprint (1999).
\bibitem{Denley} 
T. Denley,
{\it On a result of Lemke and Kleitman},
Comb., Prob. and Comput. {\bf 6} (1997), 39--43.
\bibitem{Feller} 
W. Feller,
{\it An Introduction to Probability Theory and its Applications
	(vol. I, 3rd ed.)},
Wiley, New York (1968), 31--33.
\bibitem{Foster} 
J.A. Foster and H.S. Snevily, 
{\it The 2-pebbling property and a conjecture of Graham's}, 
preprint (1995).
% \bibitem{GNN} 
% C.D. Godsil, J. Ne\v set\v ril and R. Nowakowski,
% {\it The chromatic connectivity of graphs},
% Graphs and Combin. {\bf 4} (1988), 229--233.
\bibitem{Harper} 
L.H. Harper,
{\it The morphology of partially ordered sets},
J. Combin. Theory (A) {\bf 17} (1974), 44--58.
\bibitem{HH} 
D. Herscovici and A. Higgins,
{\it The pebbling number of $C_5\Box C_5$},
Disc. Math. {\bf 187} (1998), 123--135.
\bibitem{Hsieh} 
W.N. Hsieh and D.J. Kleitman,
{\it Normalized matching in direct products of partial orders},
Stud. Appl. Math. {\bf 52} (1973), 285--289.
\bibitem{Hurl} 
G.H. Hurlbert,
{\it Two pebbling theorems},
Congress. Numer. (to appear).
\bibitem{HKT} 
G.H. Hurlbert, A.V. Kostochka and L.A. Talysheva,
{\it The dimension of interior levels of the boolean lattice},
Order {\bf 11} (1994), 29--40.
\bibitem{Kat} 
G.O.H. Katona,
{\it A theorem on finite sets}, 
in Theory of Graphs (P. Erd\H os and G.O.H. Katona, eds.),
Akad\'emiai Kiad\'o, Budapest (1968), 187--207.
\bibitem{Kr} 
J.B. Kruskal,
{\it The number of simplices in a complex}, 
in Mathematical Optimization Techniques,
Univ. California Press, Berkeley (1963), 251--278.
\bibitem{KL} 
P. Lemke and D.J. Kleitman,
{\it An addition theorem on the integers modulo $n$},
J. Number Th. {\bf 31} (1989), 335--345.
\bibitem{Lov} 
L. Lov\'asz,
{\it Combinatorial Problems and Exercises},
North Holland Pub., Amsterdam, New York, Oxford (1979).
\bibitem{Lub} 
D. Lubell,
{\it A short proof of Sperner's lemma},
J. Combin. Theory {\bf 1} (1966), 299.
\bibitem{Mac} 
F.S. Macaulay,
{\it Some properties of enumeration in the theory of modular systems},
Proc. Lond. Math. Soc. {\bf 26} (1927), 531--555.
\bibitem{Mesh} 
L.D. Meshalkin,
{\it A generalisation of Sperner's theorem on the number of subsets 
	of a finite set} (in Russian), 
Teor. Probab. Ver. Primen. {\bf 8} (1963), 219--220;
English translation in
Theory of Probab. and its Apps. {\bf 8} (1964), 204--205.
\bibitem{Moews} 
D. Moews, 
{\it Pebbling graphs},
J. Combin. Theory (B) {\bf 55} (1992), 244--252.
\bibitem{M} 
D. Moews, 
{\it personal communication} (1997).
\bibitem{Snevily} 
L. Pachter, H.S. Snevily, and B. Voxman,
{\it On pebbling graphs}, 
Congressus Numerantium {\bf 107} (1995), 65--80.
\bibitem{Saks} 
M. Saks,
{\it personal communication},
(1996).
\bibitem{Wang} 
S. Wang,
{\it Pebbling and Graham's Conjecture},
Discrete Math., {\it to appear} (1999).
\bibitem{Yam} 
K. Yamamoto,
{\it Logarithmic order of free distributive lattices},
J. Math. Soc. Japan {\bf 6} (1954), 343--353.
\end{thebibliography}
%   

% #######################################################################
% #######################################################################
%       
%       END OF PAPER
%       
\end{document}